\documentclass[11pt]{amsart}
\usepackage{listings}
\usepackage{verbatim}
\usepackage{diagbox}
\usepackage{booktabs}
\usepackage{amsmath}
\usepackage{enumerate}
\usepackage{amssymb}
\usepackage{color}
\usepackage{url}
\usepackage{graphicx}
\usepackage{fullpage}
\usepackage{tikz}
\usepackage{multirow}
\usepackage{array}
\usepackage{multirow}
\usepackage{booktabs}
\usepackage{graphicx}
\usepackage{makecell}

\newcolumntype{C}{>{\centering\arraybackslash}p{3.5cm}}

\usepackage{array, longtable}

\usepackage{amsthm,thmtools}
\theoremstyle{plain}

\newtheorem*{theorem*}{Theorem}

\usetikzlibrary{matrix}

\usepackage{xtab}
\usepackage{float}
\usepackage{color,longtable}

\usepackage[pdfencoding=auto]{hyperref}
\usepackage{bookmark}

 \theoremstyle{plain}
\newtheorem{theorem}{Theorem}[section]
\newtheorem{corollary}{Corollary}[section]

\newtheorem{lemma}{Lemma}[section]

\theoremstyle{definition}

\theoremstyle{remark}
\newtheorem{remark}{Remark}

\newtheorem*{thm*}{\textbf{Theorem}}

\numberwithin{equation}{section}

\setbox0=\hbox{$+$}
\newdimen\plusheight
\plusheight=\ht0
\def\+{\;\lower\plusheight\hbox{$+$}\;}

\setbox0=\hbox{$-$}
\newdimen\minusheight
\minusheight=\ht0
\def\-{\;\lower\minusheight\hbox{$-$}\;}

\setbox0=\hbox{$\cdots$}
\newdimen\cdotsheight
\cdotsheight=\plusheight
\def\cds{\lower\cdotsheight\hbox{$\cdots$}}

\DeclareMathOperator{\ord}{ord}
\newcommand{\Z}{\mathbb{Z}}
\newcommand{\N}{\mathbb{N}}

\makeatletter
\thmt@define@thmuse@key{tag}{%
  \thmt@suspendcounter{\thmt@envname}{#1}%
}
\makeatother

\usepackage{lipsum}

\usepackage{placeins}

\newcolumntype{C}{>{$}c<{$}}

\begin{document}

\title[Sign-patterns of Certain Infinite Products]
  {Sign-patterns of Certain Infinite Products}

\author{Zeyu Huang, Timothy Huber, James McLaughlin, Pengjun Wang, Yan Xu and Dongxi Ye}

\address{
School of Mathematics (Zhuhai), Sun Yat-sen University, Zhuhai 519082, Guangdong,
People's Republic of China}
\email{huangzy255@mail2.sysu.edu.cn}

\address{
School of Mathematical and Statistical Sciences, University of Texas Rio Grande
Valley, Edinburg, Texas 78539, USA}
\email{timothy.huber@utrgv.edu}

\address{Mathematics Department,
 25 University Avenue,
West Chester University, West Chester, PA 19383}
\email{jmclaughlin2@wcupa.edu}

\address{
School of Mathematics (Zhuhai), Sun Yat-sen University, Zhuhai 519082, Guangdong,
People's Republic of China}
\email{wangpj233@mail2.sysu.edu.cn}


\address{
School of Mathematics (Zhuhai), Sun Yat-sen University, Zhuhai 519082, Guangdong,
People's Republic of China}
\email{xuyan73@mail2.sysu.edu.cn}

\address{
School of Mathematics (Zhuhai), Sun Yat-sen University, Zhuhai 519082, Guangdong,
People's Republic of China}

\email{yedx3@mail.sysu.edu.cn}


\thanks{Dongxi Ye was supported by the Guangdong Basic and Applied Basic Research Foundation (Grant No. 2024A1515030222). Zeyu Huang and Pengjun Wang were supported by the Innovation and Enterprise Program for College Students.}

\date{\today}
\allowdisplaybreaks
\begin{abstract}
The signs of Fourier coefficients of certain eta quotients are determined by dissecting expansions for theta functions and by applying a general dissection formula for certain classes of quintuple products. A characterization is given for the coefficient sign patterns for
\[
\frac{(q^i;q^i)_{\infty}}{(q^p;q^p)_{\infty}}
\]
for integers \( i > 1 \) and primes \( p > 3 \). The sign analysis for this quotient addresses and extends a conjecture of Bringmann et al. for the coefficients of \( (q^2;q^2)_{\infty}(q^5;q^5)_{\infty}^{-1} \). The sign distribution for additional classes of eta quotients is considered. This addresses multiple conjectures posed by Bringmann et al.
  \end{abstract}

\maketitle
\allowdisplaybreaks

\section{Introduction}\label{secintro}

The topic of vanishing coefficients in the series expansion of infinite $q$-products, and the periodicity of the signs of the coefficients, goes back at least as far as the paper \cite{RS78} by
Richmond and Szekeres.

Let the sequence $\{c_n\}$ be defined by the Ramanujan product
\[
R(q):=\frac{(q^2,q^3;q^5)_{\infty}}{(q,q^4;q^5)_{\infty}} =: \sum_{n=0}^{\infty}c_n q^n,
\]
where here and subsequently the standard notation $(a_1,  \dots, a_j; q)_\infty := (a_1;q)_\infty  \cdots (a_j;q)_\infty$ is employed.
In that paper \cite{RS78}, the authors used first the saddle point method and then the circle method of Hardy and Ramanujan, as modified by Rademacher, to determine Hardy-Ramanujan-Rademacher-type series for the $c_n$, and were thus able to deduce, for $n$ sufficiently large, that
\begin{equation}\label{cnieq1}
  c_{5n}>0, \hspace{10pt}  c_{5n+1}>0, \hspace{10pt}   c_{5n+2}<0, \hspace{10pt}   c_{5n+3}<0, \hspace{10pt}   c_{5n+4}<0. 
\end{equation}
They also proved a similar result for the coefficients in the series expansion of $1/R(q)$.

Another of their results was to show that if the sequence $\{d_n\}$ is defined by 
\[
F(q):=\frac{(q^3,q^5;q^8)_{\infty}}{(q,q^7;q^8)_{\infty}} =: \sum_{n=0}^{\infty}d_n q^n,
\]
then $d_{4n+3}=0$, and that a similar vanishing coefficient result held for the coefficients in the series expansion of $1/F(q)$. 

Generalizations and extensions of those vanishing coefficient results to infinite families of infinite products were subsequently given in papers by 
Andrews and Bressoud \cite{AB79},
Alladi and Gordon  \cite{AG94} and the third author of the present paper \cite{McL15}, where the method of proof in all cases employed a special case of the Ramanujan $_1\psi_1$ summation formula.

Andrews \cite{A81} gave a combinatorial proof of the inequalities at \eqref{cnieq1}, but of greater relevance to the methods in the present paper was the proof of Hirschhorn \cite{H98}, who determined the 5-dissection of $R(q)$, i.e. expressed $R(q)$ in the form
\[
R(q)= R_0(q^5)+q R_1(q^5)-q^7 R_2(q^5)-q^3 R_5(q^5)-q^{14} R_4(q^5),
\]
where the coefficients in the series expansion of each of the $R_i(q^5)$, $0 \leq i \leq 4$, are all nonnegative and strictly positive eventually, making the sign periodicity stated at \eqref{cnieq1} obvious.

A new class of  infinite $q$-products with the property that when the product is expanded as a series in $q$, then the coefficients in one or more arithmetic progressions vanish, was introduced by Hirschhorn in  \cite{H19}, where he showed that if the sequences $\{a_n\}$ and $\{b_n\}$ are defined by
 {\allowdisplaybreaks
\begin{align*}
\sum_{n=0}^{\infty}
a_n q^n &:= (-q,-q^4 ;q^5 )_{\infty} (q,q^9 ;q^{10})_{\infty}^3,\\
\sum_{n=0}^{\infty}
b_n q^n &:= (-q^2,-q^3 ;q^5 )_{\infty} (q^3,q^7 ;q^{10})_{\infty}^3.
\end{align*}}
Then $a_{5n+2}=a_{5n+4}=b_{5n+1}=b_{5n+4}=0$.  Similar results were subsequently  proven by Tang \cite{T19},  Baruah and Kaur   \cite{BK19}, the third author of the present paper \cite{M19a} and with Zimmer \cite{McLZ22}. Motivated by the results and methods of  Richmond and Szekeres  \cite{RS78} and Hirschhorn \cite{H19} and others,  the topic of vanishing coefficients in the series expansion of various kinds of infinite $q$-products is an active area of research - see for example the papers of  
Channabasavayya and Dasappa \cite{CD24},  
Channabasavayya,  Keerthana and Dasappa \cite{CKD24},
Chern and Tang \cite{CT20, CT21, CT24}, 
Daniels \cite{D24},
Dou and Xiao \cite{DX21},   
Kaur and Vanda \cite{KV22, VK22}, 
Liu \cite{L24}, 
Rajkhowa and Saikia \cite{RS24},
Somashekara and Thulasi \cite{ST23} and 
    Tang \cite{T19b, T23, T23b, T24}.
    
    The topic of periodicity of the signs of the coefficients in the series expansion of infinite products has seen less investigation, but a number of authors have obtained results by deriving $m$-dissections of infinite products (for various integers $m>1$) as Hirschhorn \cite{H98} did - see for example the papers of 
  Chern and Tang \cite{CT21}, 
    Dou and Xiao \cite{DX21},  
  Tang   \cite{T21} and 
Xia and Zhao \cite{XZ22}.

Another important result was that of  Andrews \cite{GE1}, where he proved the signs of the Fourier coefficients of the infinite Borwein product
\begin{equation}\label{GPeq1}
    G_{p} (q) =\frac{(q;q)_{\infty}}{(q^{p}; q^{p})_{\infty}} = \sum_{n=0}^{\infty} c_{p}(n)q^{n}
\end{equation}
 are periodic modulo $p$. Andrews also remarked that Frank Garvan and Peter Borwein had a different (unpublished) proof of this sign periodicity of $G_p$.

The function $G_p$ was further studied in a recent paper of 
Schlosser and Zhou \cite{SZ23}, where they made some conjectures about the periodicity of the signs in the series expansion of $G_p^{\delta}$ for real $\delta$ in certain intervals, and proved, for $0.227 \leq \delta \leq 2.9999$ and $n \geq 158$, that the signs of the coefficients $c_3^{\delta}$ in the series expansion of $G_3^{\delta}$  have period 3. 

For reasons of compactness, write $(q^{j}; q^{j})_{\infty} =: f_j$, so that $G_p = f_1 f_p^{-1}$. The study of $G_p$ has led to the investigation of vanishing coefficients and periodicity of sign changes in other eta quotients 
\[
\prod_{j=1}^{m}f_j^{\delta_j}=:\sum_{n\ge0}C_{1^{\delta_1} 2^{\delta_2}\cdots m^{\delta_m}}(n)q^n, \, \delta_j \in \Z,
\]
for example in the recent papers by Bringmann et. al. \cite{Br,BHHK}, where we have employed their notation for the coefficients. If we continue with their notation and define the \begin{it}vanishing set\end{it}
\[
\mathcal{S}_{1^{\delta_1} 2^{\delta_2}\cdots m^{\delta_m}}:=\left\{n\in\N: C_{1^{\delta_1} 2^{\delta_2}\cdots m^{\delta_m}}(n)=0\right\},
\]
then an example of one of their results from \cite{BHHK2} is that 
\[
\mathcal{S}_{1^{-1}2^{10}3^{-1}4^{-4}}=\mathcal{S}_{1^72^{-2}3^{-1}}=\left\{n\in\N: n\equiv 2\pmod{3}\text{ and }\exists p\equiv 3\pmod{4}, \ord_p(n)\text{ {\rm is odd}}\right\}.
\]
An example of one of their results on sign periodicity from \cite{BHHK} is the following:\\
The signs of the sequence $\{C_{1^{4}2^{2}4^{-2}}(n)\}_{n\ge1}$ have period $8$. In particular, 
\[
	C_{1^{4}2^{2}4^{-2}}(n)
	\begin{cases}
		>0&\text{if }n\equiv 0,3,7\pmod{8},\\
		<0&\text{if }n\equiv 1,4,5\pmod{8},\\
		=0&\text{if }n\equiv 2\pmod{4}.
	\end{cases}
\]
In the paper \cite{BHHK3} the authors prove some additional results, including the following for period 9:
\[
		C_{1^93^{-5}}(n)=\begin{cases}
			>0&\text{if }n\equiv 0,2,5,6,8\pmod{9},\\
			<0&\text{if }n\equiv 1,3,4,7\pmod{9}
		\end{cases}
		\]
and make a large number of conjectures about the sign periodicity of the coefficients in the series expansion of particular eta quotients; for example, they conjectured \cite[Table 1]{BHHK3} that
\begin{equation}
    \label{C25}
     C_{2^{1}5^{-1}}\begin{cases}
        >0&\mbox{if $n\equiv0\pmod{5}$,}\\
        <0&\mbox{if $n\equiv 2,4\pmod{5}$,}\\
        =0&\mbox{otherwise.}
    \end{cases}
\end{equation}


Inspired by Bringmann et al.'s works and their conjectures, one of the primary goals of this paper is to prove the coefficient sign pattern claimed in the next theorem.






\begin{theorem}\label{main1}
    Let $p>3$ be a prime, and let $i>1$ be an integer not divisible by~$p$. Write
    \begin{align}\label{neqremove1} 
       \frac{(q^i;q^i)_{\infty}}{(q^p;q^p)_{\infty}}=\sum_{n=0}^{\infty}a_n q^n.
    \end{align}
    For $0\leq r\leq p-1$, when $p\equiv1\pmod{3}$, let
    $$
   \mathcal{L}(r)=
  \begin{cases}
      6r^2+r    & \mbox{if } r\leq \frac{4p-1}{12}, \\
      6r^2+r-8pr+\frac{8p^2-2p}{3} & \mbox{if } \frac{4p-1}{12}<r\leq \frac{10p-1}{12}, \\
      6r^2+r-12pr+(6p^2-p) & \mbox{if } r>\frac{10p-1}{12},
  \end{cases}
    $$
    and when $p\equiv-1\pmod{3}$, let
    $$
     \mathcal{L}(r)=
  \begin{cases}
      6r^2+r    & \mbox{if } r\leq \frac{2p-1}{12}, \\
      6r^2+r-4pr+\frac{2p^2-p}{3} & \mbox{if } \frac{2p-1}{12}<r\leq \frac{8p-1}{12}, \\
      6r^2+r-12pr+(6p^2-p) & \mbox{if } r>\frac{8p-1}{12}.
  \end{cases}
    $$
    Define 
    $$
    N=N(p,i)=\max\left(\bigcup_{s=0}^{p-1}\min(i\mathcal{L}(r):\,i\mathcal{L}(r)\equiv s\pmod{p},\,0\leq r\leq p-1)\right)-p.
    $$
  Then for $n>N$, 
  \begin{enumerate}
      \item  if $p \equiv 1 \pmod 3 $, 
   \begin{align*}
       a_n\begin{cases}
           >0, &\mbox{if $n\equiv i(6r^2+r) \pmod p$ with $r\leq \frac{4p-1}{12}$ or $r>\frac{10p-1}{12}$,}\\
           <0, &\mbox{if $n\equiv i(6r^2+r) \pmod p$ with $\frac{4p-1}{12} <r\leq \frac{10p-1}{12}$,}\\
           =0, &\mbox{if  $n\neq i(6r^2+r) \pmod p$,}
       \end{cases}
   \end{align*}
   \item 
    if $p \equiv -1 \pmod 3 $, then
   \begin{align*}
       a_n\begin{cases}
           >0, &\mbox{if $n\equiv i(6r^2+r) \pmod p$ with $r\leq \frac{2p-1}{12}$ or $r>\frac{8p-1}{12}$,}\\
           <0, &\mbox{if $n\equiv i(6r^2+r) \pmod p$ with $\frac{2p-1}{12} <r\leq \frac{8p-1}{12}$,}\\
           =0, &\mbox{if  $n\neq i(6r^2+r) \pmod p$.}
       \end{cases}
   \end{align*}
  \end{enumerate}

\end{theorem}

\begin{remark}
   For the case of $i$ being divisible by~$p$, it is straightforward to show that
   $$
   a_n\begin{cases}
        >0 &\mbox{if $p|n$,}\\
        =0 &\mbox{otherwise.}
    \end{cases}
    $$
    We leave the details to the reader.
\end{remark}

Specializing $p=5$ and $i=2$ in Theorem~\ref{main1} yields~\eqref{C25} and thus, 
\begin{corollary}
      Write
    $$
    \frac{(q^{2};q^{2})_{\infty}}{(q^{5};q^{5})_{\infty}}=\sum_{n=0}^{\infty}a_{n}q^{n}.
    $$
    Then
    $$
    a_{n}\begin{cases}
        >0&\mbox{if $n\equiv0\pmod{5}$,}\\
        <0&\mbox{if $n\equiv 2,4\pmod{5}$,}\\
        =0&\mbox{otherwise.}
    \end{cases}
    $$
\end{corollary}

In general, as is indicated in Theorem~\ref{main1}, the sign-periodicity of the coefficients $a_{n}$ of $\frac{(q^i;q^i)_{\infty}}{(q^p;q^p)_{\infty}}$ is not universal for any~$n$. However, the lower bound~$N$ on $n$ given in the theorem is a sharp one. In the following table, we display all the examples for $2\leq i\leq4$ and $5\leq p\leq 13$ that follow from the theorem. In a cell corresponding to a pair $(p,i)$, the upper row represents the sign-periodicity of $a_{n}$ in order with $n\equiv 0,\ldots, p-1\pmod{p}$, and the lower row gives the largest subset of $n$'s for which the sign-periodicity holds.
\begin{center}
    \begin{tabular}{c|C|C|C}
\toprule
\diagbox[width=3em, height=2em]{$p$}{$i$} & \textbf{2} & \textbf{3} & \textbf{4} \\
\midrule
\multirow{2}{*}{\textbf{5}} 
  & \makecell{+0-0-} 
  & \makecell{+-0-0} 
  & \makecell{+00--} \\
  & \makecell{n \geq 0} 
  & \makecell{n \geq 2} 
  & \makecell{n \geq 4} \\
\midrule
\multirow{2}{*}{\textbf{7}} 
  & \makecell{+0-+-00} 
  & \makecell{++0-00-} 
  & \makecell{+-00-0+} \\
  & \makecell{n \geq 4} 
  & \makecell{n \geq 9} 
  & \makecell{n \geq 14} \\
\midrule
\multirow{2}{*}{\textbf{11}} 
  & \makecell{+0-+-000-0+} 
  & \makecell{+-0-+0-000+} 
  & \makecell{+000--+0-+0} \\
  & \makecell{n \geq 20} 
  & \makecell{n \geq 35} 
  & \makecell{n \geq 50} \\
\midrule
\multirow{2}{*}{\textbf{13}} 
  & \makecell{++-0-+0000+-0} 
  & \makecell{+++-00-0+0-00} 
  & \makecell{+0+0-00+--+00} \\
  & \makecell{n \geq 32} 
  & \makecell{n \geq 54} 
  & \makecell{n \geq 76} \\
\bottomrule
\end{tabular}
\end{center}

In addition to Theorem~\ref{main1}, we also formulate the sign-periodicity of the coefficients of multiple particular infinite products in the following theorem.

\begin{theorem}\label{main2}
    \begin{enumerate}

     \item  For $p\geq3$ a prime,  write
    $$
    \frac{(q^2;q^2)^2_{\infty}}
    {(q;q)_{\infty}(q^{p};q^{p})_{\infty}}=\sum_{n=0}^{\infty}a_{n}q^{n}.
    $$
    Let
    $$
    N=N(p)=\max\left(\bigcup_{s=0}^{p-1}\min\left(\frac{r(r+1)}{2}:\,\frac{r(r+1)}{2}\equiv s\pmod{p},\,0\leq r\leq p-1\right)\right)-p.
    $$
    Then for $n>N$,
    $$
    a_{n}\begin{cases}
        >0&\mbox{if $n\equiv \frac{r(r+1)}{2}\pmod{p}$ for some $r\in\mathbb{Z}/p\mathbb{Z}$,}\\
        =0&\mbox{otherwise.}
    \end{cases}
    $$

\item  For $p=1$ or an odd prime,  write 
    $$
    \frac{(q;q)^2_{\infty}}{(q^2;q^2)_\infty(q^{4p};q^{4p})_{\infty}}=\sum_{n=0}^{\infty}a_{n}q^{n}.
    $$
    Let
    $$
    N=N(p)=\max\left(\bigcup_{s=0}^{4p-1}\min\left(r^{2}:\,r^{2}\equiv s\pmod{p},\,0\leq r\leq 4p-1\right)\right)-4p.
    $$
    Then for $n>N$,
    $$
    a_{n}\begin{cases}
        <0&\mbox{if $n\equiv 4t^{2}+4t+1\pmod{4p}$ for some $t\in\mathbb{Z}/p\mathbb{Z}$,}\\
        >0&\mbox{if $n\equiv 4t^{2}\pmod{4p}$  for some $t\in\mathbb{Z}/p\mathbb{Z}$,}\\
        =0&\mbox{otherwise, i.e.,  $n\equiv 2,3\pmod{4}$ or $\left(\frac{n}{p}\right)=-1$.}
        \end{cases}
    $$

\item    Write
    $$
    \frac{(q;q)^3_{\infty}}{(q^3;q^3)_{\infty}^2}=\sum_{n=0}^{\infty}a_{n}q^{n}.
    $$
    Then
    $$
    a_{n}\begin{cases}
        >0&\mbox{if $n\equiv 0\pmod{3}$,}\\
        <0&\mbox{if $n\equiv 1\pmod{3}$,}\\
        =0&\mbox{otherwise.}
    \end{cases}
    $$

    \item Write
$$
\frac{(q;q)^2_{\infty}}
{(q^2;q^2)_\infty(q^3;q^3)_{\infty}^{2}} = \sum_{n=0}^{\infty}a_{n}q^{n}
$$
Then
$$
a_{n}\begin{cases}
        >0&\mbox{if $n\equiv 0\pmod{3}$,}\\
        <0&\mbox{if $n\equiv 1\pmod{3}$,}\\
        =0&\mbox{otherwise.}
    \end{cases}
$$

  \item      Write
    $$
    \frac{(q;q)^4_{\infty}}{(q^2;q^2)^2_\infty(q^4;q^4)_{\infty}}=\sum_{n=0}^{\infty}a_{n}q^{n}.
    $$
    Then
    $$
    a_{n}\begin{cases}
        >0&\mbox{if $n\equiv 0,2\pmod{4}$,}\\
        <0&\mbox{if $n\equiv 1\pmod{4}$,}\\
        =0&\mbox{otherwise.}
        \end{cases}
    $$

\item 

Write
$$
\frac{(q^2;q^2)^{10}_{\infty}}{(q;q)_{\infty}^4(q^4;q^4)_{\infty}^5}
=\sum_{n=0}^{\infty}a_{n}q^{n}.
$$
Then
    $$
    a_{n}\begin{cases}
        >0&\mbox{if $n\equiv 0,1,2\pmod{4}$,}\\
        =0&\mbox{otherwise.}
    \end{cases}
    $$


           \item  Write
    $$
    \frac{(q;q)^2_{\infty}}
    {(q^{5};q^{5})^3_{\infty}}=\sum_{n=0}^{\infty}a_{n}q^{n}.
    $$
    Then
    $$
    a_{n}\begin{cases}
        >0&\mbox{if $n\equiv 0,3,4\pmod{5}$,}\\
        <0&\mbox{otherwise.}
    \end{cases}
    $$

\item  
Write
\begin{align*}
        \frac{(q;q)^9}{(q^3;q^3)^9}=\sum_{n=0}^\infty a_nq^n.
    \end{align*}
    Then
    $$
a_{n}\begin{cases}
    >0&\mbox{if $n\equiv 0,2,5,8\pmod{9}$,}\\
    <0&\mbox{if $n\equiv 1,3,4,7\pmod{9}$,}\\
    =0&\mbox{otherwise.}
\end{cases}
    $$

\item Let $i\geq11$ be an integer, and write
    \begin{align*}
        \frac{(q;q)^9}{(q^3;q^3)^i}=\sum_{n=0}^\infty a_nq^n.
    \end{align*}
Then
\begin{enumerate}
    \item when $i=11$, 
    $$
    a_{n}\begin{cases}
    >0&\mbox{if $n\equiv 0,2,5,8\pmod{9}$,}\\
    <0&\mbox{otherwise,}
    \end{cases}
    $$ 

    \item when $i=12$, 
    $$
    a_{n}\begin{cases}
    >0&\mbox{if $n\equiv 0,2,5,8\pmod{9}$,}\\
    <0&\mbox{if $n\equiv 1,4,7\pmod{9}$,}\\
    =0&\mbox{otherwise,}
    \end{cases}
    $$ 

    \item when $i>12$,
    $$
    a_{n}\begin{cases}
    >0&\mbox{if $n\equiv 0,2\pmod{3}$,}\\
    <0&\mbox{otherwise.}
    \end{cases}
    $$ 
\end{enumerate}
    
    
    
    
    \end{enumerate}
\end{theorem}

\begin{remark}
A number of particular conjectures posed in \cite[Table 1]{BHHK3} are covered by Theorem~\ref{main2}:
\begin{itemize}
    \item When $p=5$, part (1) together with checking the first $\frac{p(p+1)}{2}=15$ terms proves the case corresponding to {$5/++0+0$}.

    \item   When $p=1$, part (2) proves the case corresponding to {$4/+-00$}.

    \item Parts (3) and (4) respectively prove the first- and second-to-the-right case corresponding to {$3/+0-$}.


    \item Part (5) proves the case corresponding to {$4/+-+0$}.

    \item Part (6) proves the upper case in the first cell corresponding to {$4/+++0$}.

    \item Part (7) proves the middle case in the second cell corresponding to $5/+--++$.

    \item Part (8) proves the case corresponding to $9/+-+--+0-+$.
\end{itemize}

\end{remark}

\section{Nuts and Bolts}

In this section, we shall state and prove some preliminary results that will be useful in proving Theorems~\ref{main1} and~\ref{main2}. We start with the following technical lemma that will be used frequently throughout the work. The proof of the lemma is straightforward and left to the reader.

\begin{lemma}\label{technical}
\begin{enumerate}
    \item    For $|q|<1$, and any two subsets $S_{1},S_{2}$ of $\mathbb{Z}_{>0}$ such that $S_{1}\subset S_{2}$, one has
    $$
    \frac{\prod_{n\in S_{1}}(1-q^{n})}{\prod_{n\in S_{2}}(1-q^{n})}=\frac{1}{\prod_{n\in S_{2}-S_{1}}(1-q^{n})}.
    $$
    Moreover, if one writes
    $$
    \frac{1}{\prod_{n\in S_{2}-S_{1}}(1-q^{n})}=1+\sum_{n=1}^{\infty}a_{n}q^{n},
    $$
    then $a_{n}\geq0$ for any~$n\geq1$.

    \item For an integer $m\geq2$ and any series $1+\sum_{n=1}^{\infty}b_{n}q^{mn}$ with $b_{n}\geq0$, write
    $$
    \left(1+\sum_{n=1}^{\infty}b_{n}q^{mn}\right)\frac{1}{1-q^{m}}=\sum_{n=0}^{\infty}a_{n}q^{mn}.
    $$
    Then $a_{n}\geq1$ for any~$n\geq0$.
\end{enumerate}
 
\end{lemma}

The remainder of the present section is split into three subsections in accordance with the nature of the preliminary results to be displayed.

\subsection{Dissection formulas}

The following $m$-dissection formula can be found in~\cite{HMY}.

\begin{theorem}\label{qqqL}
For $M\geq3$, $1\leq j<M/2$, and $m\equiv\pm1\pmod{3}$, one has that
   \begin{align}\label{neqremove2}
   &(q^{j},q^{M-j},q^{M};q^{M})_{\infty}(q^{M-2j},q^{M+2j};q^{2M})_{\infty}\\
   &=\sum_{r=0}^{m-1}(-1)^{s(r)} q^{\mathcal{L}(r)}(q^{t_{1}(r)},q^{m^{2}M-t_{1}(r)},q^{m^{2}M};q^{m^{2}M})_{\infty}(q^{t_{2}(r)},q^{2m^{2}M-t_{2}(r)};q^{2m^{2}M})_{\infty},\nonumber
 \end{align}
 where
 \begin{align*}
  t_{1}(r)&=m M (m\pm(6 r-1))/6\pm j m\pmod{m^{2}M},\\
    t_{2}(r)&=m^2M +2 j m\pm M (m\pm(6 r-1)) m/3\pmod{2m^{2}M},\\
     \mathcal{L}(r)&=\frac{7m^{2}M}{24}+\frac{1}{2}t_{1}(r)\left(\frac{t_{1}(r)}{m^{2}M}-1\right)+\frac{1}{2}t_{2}(r)\left(\frac{t_{2}(r)}{2m^{2}M}-1\right)\\
     &\quad-\left(\frac{7M}{24}+\frac{1}{2}j\left(\frac{j}{M}-1\right)+\frac{1}{2}(M-2j)\left(\frac{M-2j}{2M}-1\right)\right),
 \end{align*}
 and if $m\equiv 1 \pmod 3$, then
\begin{equation}\label{seq1}
 s(r)=
 \begin{cases}
   0, & \mbox{if } r\leq\frac{(2 m+1) M-6 j}{6 M}, \\
  1, & \mbox{if } \frac{(2 m+1) M-6 j}{6 M}<r\leq\frac{(5 m+1) M-6 j}{6 M}, \\
  2, & \mbox{if } r> \frac{(5 m+1) M-6 j}{6M},
 \end{cases}
\end{equation}
and if $m\equiv -1 \pmod 3$, then
\begin{equation}\label{seq2}
 s(r)=
 \begin{cases}
   0, & \mbox{if } r\leq \frac{(m+1) M-6 j}{6 M}, \\
  1, & \mbox{if } \frac{(m+1) M-6 j}{6 M}<r\leq \frac{(4 m+1) M-6 j}{6 M}, \\
  2, & \mbox{if } r>\frac{(4 m+1) M-6 j}{6  M}.
 \end{cases}
\end{equation}

\end{theorem}

An immediate implication of Theorem~\ref{qqqL} is the following explicit dissection formula for $(q;q)_{\infty}$.

\begin{corollary}
\label{qqdissec}
 For  an $m\equiv\pm1\pmod{3}$, one has that  
  \begin{align}\label{neqremove3}
   &(q;q)_{\infty}=\sum_{r=0}^{m-1}(-1)^{s(r)} q^{\mathcal{L}(r)}(q^{t_{1}(r)},q^{4m^{2}-t_{1}(r)},q^{4m^{2}};q^{4m^{2}})_{\infty}(q^{t_{2}(r)},q^{8m^{2}-t_{2}(r)};q^{8m^{2}})_{\infty},\nonumber
 \end{align}
where if $m\equiv 1 \pmod 3$, 
\begin{align*}
  t_{1}(r)&=
  \begin{cases}
      \frac{2m^2+m}{3}+4mr, & \mbox{if } r<\frac{10m-1}{12}, \\
      \frac{-10m^2+m}{3}+4mr, & \mbox{if } r>\frac{10m-1}{12},
  \end{cases} \\
  t_{2}(r)&=
  \begin{cases}
      \frac{16m^2+2m}{3}+8mr, & \mbox{if } r<\frac{4m-1}{12}, \\
      \frac{-8m^2+2m}{3}+8mr, & \mbox{if } r>\frac{4m-1}{12},
  \end{cases} \\
  \mathcal{L}(r)&=
  \begin{cases}
      6r^2+r    & \mbox{if } r\leq \frac{4m-1}{12}, \\
      6r^2+r-8mr+\frac{8m^2-2m}{3} & \mbox{if } \frac{4m-1}{12}<r\leq \frac{10m-1}{12}, \\
      6r^2+r-12mr+(6m^2-m) & \mbox{if } r>\frac{10m-1}{12},
  \end{cases}\\
   s(r)&=
 \begin{cases}
   0, & \mbox{if } r\leq \frac{4m-1}{12}, \\
  1, & \mbox{if } \frac{4m-1}{12}<r\leq \frac{10m-1}{12}, \\
  2, & \mbox{if } r>\frac{10m-1}{12},
 \end{cases}
\end{align*} 
 and  if $m\equiv -1 \pmod 3$, 
\begin{align*}
  t_{1}(r)&=
  \begin{cases}
      \frac{2m^2-m}{3}-4mr, & \mbox{if } r<\frac{2m-1}{12}, \\
      \frac{14m^2-m}{3}-4mr, & \mbox{if } r>\frac{2m-1}{12},
  \end{cases} \\
  t_{2}(r)&=
  \begin{cases}
      \frac{8m^2+2m}{3}+8mr, & \mbox{if } r<\frac{8m-1}{12}, \\
      \frac{-16m^2+2m}{3}+8mr, & \mbox{if } r>\frac{8m-1}{12},
  \end{cases} \\
  \mathcal{L}(r)&=
  \begin{cases}
      6r^2+r    & \mbox{if } r\leq \frac{2m-1}{12}, \\
      6r^2+r-4mr+\frac{2m^2-m}{3} & \mbox{if } \frac{2m-1}{12}<r\leq \frac{8m-1}{12}, \\
      6r^2+r-12mr+(6m^2-m) & \mbox{if } r>\frac{8m-1}{12},
  \end{cases}\\
   s(r)&=
 \begin{cases}
   0, & \mbox{if } r\leq \frac{2m-1}{12}, \\
  1, & \mbox{if } \frac{2m-1}{12}<r\leq \frac{8m-1}{12}, \\
  2, & \mbox{if } r>\frac{8m-1}{12}.
 \end{cases}
\end{align*}

\end{corollary}

\begin{proof}
    This follows from taking $M=4$, $ j=1$ in Theorem~\ref{qqqL} and noticing that
    $$
    (q;q)_{\infty}=(q^1,q^3,q^4;q^4)_{\infty}(q^2,q^6;q^8)_{\infty}.
    $$
In this case,    when $m\equiv1\pmod3$
    \begin{align*}
  t_{1}(r)&\equiv\frac{2m^2+m}{3}+4mr\pmod{4m^{2}},\\
    t_{2}(r)&\equiv\frac{16m^2+2m}{3}+8mr\pmod{8m^{2}}.\\
 \end{align*}
 One finds that $t_1(\frac{10m-1}{12})=4m^2<t_1(m-1)<8m^2 \text{ and } t_2(\frac{4m-1}{12})=8m^2<t_2(m-1)<16m^2$.
 So,
 \begin{align*}
 t_{1}(r)&=
  \begin{cases}
      \frac{2m^2+m}{3}+4mr, & \mbox{if } r<\frac{10m-1}{12}, \\
      \frac{-10m^2+m}{3}+4mr, & \mbox{if } r>\frac{10m-1}{12},
  \end{cases} \\
  t_{2}(r)&=
  \begin{cases}
      \frac{16m^2+2m}{3}+8mr, & \mbox{if } r<\frac{4m-1}{12}, \\
      \frac{-8m^2+2m}{3}+8mr, & \mbox{if } r>\frac{4m-1}{12}.
  \end{cases} 
 \end{align*} 
 Also,  when $m\equiv -1 \pmod 3$, 
\begin{align*}
  t_{1}(r)&\equiv\frac{2m^2-m}{3}-4mr\pmod{4m^{2}},\\
    t_{2}(r)&\equiv\frac{8m^2+2m}{3}+8mr\pmod{8m^{2}},\\
 \end{align*}
Similarly, $-4m^2<t_1(m-1)<t_1(\frac{2m-1}{12})=0<t_1(1)<4m^2 \text{ and } t_2(\frac{8m-1}{12})=8m^2<t_2(m-1)<16m^2$.
 Therefore,
\begin{align*}
  t_{1}(r)&=
  \begin{cases}
      \frac{2m^2-m}{3}-4mr, & \mbox{if } r<\frac{2m-1}{12}, \\
      \frac{14m^2+m}{3}-4mr, & \mbox{if } r>\frac{2m-1}{12},
  \end{cases} \\
  t_{2}(r)&=
  \begin{cases}
      \frac{8m^2+2m}{3}+8mr, & \mbox{if } r<\frac{8m-1}{12}, \\
      \frac{-16m^2+2m}{3}+8mr, & \mbox{if } r>\frac{8m-1}{12}.
  \end{cases} 
\end{align*}
As a consequence, $\mathcal{L}(r)$ may be formulated as stated.
\end{proof}

Note that the $m$-dissection formulas given above hold only for $m$ coprime to~$3$. We conclude the present subsection with the next theorem concerning the $3$-dissection of the infinite products~$(q;q)_{\infty}$ and $(q;q)_{\infty}^{3}$.
\begin{theorem}\label{3dissec}
The following identities hold.
    \begin{align*}
   (q;q)_{\infty}&= \frac{(q^{3};q^{3})_{\infty}}{(q^3, q^6, q^9, q^{18}, q^{21}, q^{24}, q^{27}; q^{27})_\infty}- q\frac{(q^{3};q^{3})_{\infty}}{(q^3, q^9, q^{12}, q^{15}, q^{18}, q^{24}, q^{27}; q^{27})_\infty}\\
    &\quad - q^{2}\frac{(q^{3};q^{3})_{\infty}}{(q^6, q^9, q^{12}, q^{15}, q^{18}, q^{21}, q^{27}; q^{27})_\infty},\\
         (q; q)_\infty^3 &= (q^3;q^3)_\infty\left(1+6\sum_{n=1}^{\infty}q^{3n}\frac{1-q^{3n}}{1-q^{9n}}\right)-3q (q^9; q^9)_\infty^3.
    \end{align*}

\end{theorem}

\begin{proof}
    Take $m=3$ in \cite[(2.1) and Theorem~4.1]{M}, respectively.
\end{proof}

\subsection{The Borwein cubic theta functions}

Define
\begin{align*}
    a(q)&:=\sum_{m,n=-\infty}^{\infty}q^{m^{2}+mn+n^{2}}, \quad
     b(q):=\sum_{m,n=-\infty}^{\infty}e^{2\pi i(m-n)/3}q^{m^{2}+mn+n^{2}}\,,\\
     c(q)&:=\sum_{m,n=-\infty}^{\infty}q^{(m+\frac{1}{3})^{2}+(m+\frac{1}{3})(n+\frac{1}{3})+(n+\frac{1}{3})^{2}}.
\end{align*}
These functions are called the Borwein cubic theta functions. 

Clearly the coefficients of $a(q)$ are all non-negative. Moreover, the functions $a(q),b(q), c(q)$ satisfy the following identities \cite[Chapter 3]{C}:
\begin{equation}\label{abc1}
     b(q)=\frac{(q;q)_{\infty}^{3}}{(q^{3};q^{3})_{\infty}}=a(q^{3})-c(q^3), \qquad c(q) = 3q^{\frac{1}{3}}\frac{(q^{3};q^{3})_{\infty}^{3}}{(q;q)_{\infty}},  \qquad a^{3}(q) = b^{3}(q) + c^{3}(q).
\end{equation}
We shall make use of \eqref{abc1} to establish the following lemma that will be useful in proving items (8) and (9) of Theorem~\ref{main2}.







\begin{lemma}
\label{lem9}
    Let $i>3$ be an integer, and write
    \begin{align*}
        \frac{(q;q)_{\infty}^9}{(q^3;q^3)_{\infty}^i}=\sum_{n=0}^\infty a_nq^n.
    \end{align*}
    Then
    $$
    a_{n}\begin{cases}
        >0&\mbox{if $n\equiv2\pmod{3}$,}\\
        <0&\mbox{if $n\equiv1\pmod{3}$.}
    \end{cases}
    $$
\end{lemma}
\begin{proof}
Note by \cite[Theorem 3.15]{C} that
$$
\frac{(q^3;q^3)^3_\infty}{(q;q)_\infty}=\sum_{m_1 + m_2 + m_3 = 0} q^{3(m_1^2 + m_2^2 + m_3^2)/2 + m_1 + 2m_2 + 3m_3},
$$
so the coefficients of $\frac{(q^3;q^3)^3_\infty}{(q;q)_\infty}$ are all non-negative. Then if one writes for any positive integer $j$,
\begin{align*}
    \left(\frac{(q^3;q^3)^3_\infty}{(q;q)_\infty}\right)^{j}\frac{1}{(q^3;q^3)_\infty}=\sum_{n=0}^\infty c_nq^{3n},
\end{align*}
by Lemma~\ref{technical} one can find that $c_{n}\geq1$. From~\eqref{abc1},
\begin{align}
    \frac{(q;q)^9}{(q^3;q^3)^i}&=\frac{1}{(q^3;q^3)^{i-3}}\left(\frac{(q;q)_{\infty}^{3}}{(q^{3};q^{3})_{\infty}}\right)^3 \nonumber\\
    &=\frac{1}{(q^3;q^3)^{i-3}}\left(a(q^{3})-3q\frac{(q^{9};q^{9})_{\infty}^{3}}{(q^{3};q^{3})_{\infty}}\right)^3 \nonumber \\
    &\hspace{-15pt}=\frac{1}{(q^3;q^3)^{i-3}}\left(a^3(q^3)-9qa^2(q^3)\frac{(q^{9};q^{9})_{\infty}^{3}}{(q^{3};q^{3})_{\infty}}+27q^{2}a(q^3)\left(\frac{(q^{9};q^{9})_{\infty}^{3}}{(q^{3};q^{3})_{\infty}}\right)^2-27q^3\left(\frac{(q^{9};q^{9})_{\infty}^{3}}{(q^{3};q^{3})_{\infty}}\right)^3\right), \label{to}
\end{align}
and thus,
$$
\sum_{n=0}^{\infty}a_{3n+1}q^{n}=-a^2(q)\frac{(q^{3};q^{3})_{\infty}^{3}}{(q;q)_{\infty}}\times\frac{1}{(q^{3};q^{3})_{\infty}^{i-3}}\quad\mbox{and}\quad \sum_{n=0}^{\infty}a_{3n+2}q^{n}=27a(q)\left(\frac{(q^{3};q^{3})_{\infty}^{3}}{(q;q)_{\infty}}\right)^2\frac{1}{(q^{3};q^{3})_{\infty}^{i-3}}.
$$
The desired conclusion follows from these dissection formulas together with \eqref{to} and the non-negativity of the coefficients of $a(q)$.
\end{proof}

\subsection{Infinite product as a theta series}



As the dissection formula given in Theorem~\ref{qqqL} might not be applicable to the cases considered in Theorem~\ref{main2},
one of the tricks that we will employ in the treatment of those cases is to factor out an appropriate component that can be expressed as a theta series and $m$-dissect the theta series by analyzing its associated lattice points. This can lead us to identify the sign-periodicity of the coefficients we are concerned with. For use in our arguments, we list all the identities linking an infinite product to a theta series.

\begin{theorem}
The following identities hold.
\begin{align}
    \frac{(q;q)_{\infty}^{2}}{(q^{2};q^{2})_{\infty}} &= \sum_{n=-\infty}^{\infty} (-1)^n q^{n^2},\label{eta1}\\
    \frac{(q^2;q^2)_{\infty}^{2}}{(q;q)_{\infty}} &= \sum_{n=0}^{\infty} q^{n(n+1)/2},\label{eta2}\\
    \frac{(q^2;q^2)^5_{\infty}}{(q;q)^2_{\infty}(q^4;q^4)^2_\infty} &=\sum_{n=-\infty}^{\infty}q^{n^2},\label{eta5}\\
\frac{(q;q)_{\infty}^{2}(q^{6};q^{6})_{\infty}}{(q^{2};q^{2})_{\infty}(q^{3};q^{3})_{\infty}} &= \frac{1}{2}\sum_{n=-\infty}^{\infty}a(n)q^{\frac{n^{2}-1}{8}}, \label{eta6}
\end{align}
where $a(n)=0,1,0,-2,0,1$ in order with $n\equiv 0,1,2,3,4,5\pmod{6}$.
\end{theorem}

\begin{proof}
    See, e.g., \cite[Theorem~5.9.4]{CS} and the references wherein.
\end{proof}

\section{Proof of Theorem~\ref{main1}}
This section is devoted to proving Theorem~\ref{main1}. As both cases considered in the theorem are parallel to one another, we only detail the proof of the case of $p\equiv1\pmod{3}$ whose conclusion is now restated in Theorem~\ref{main3} for the convenience of the reader.

\begin{theorem}
\label{main3}
    Let $p\equiv1\pmod{3}$ be a prime, and let $i>1$ be an integer not divisible by~$p$. Write
    \begin{align}\label{neqremove4} 
       \frac{(q^i;q^i)_{\infty}}{(q^p;q^p)_{\infty}}=\sum_{n=0}^{\infty}a_n q^n.
    \end{align}
  For $0\leq r\leq p-1$,  let
    $$
   \mathcal{L}(r)=
  \begin{cases}
      6r^2+r    & \mbox{if } r\leq \frac{4p-1}{12}, \\
      6r^2+r-8pr+\frac{8p^2-2p}{3} & \mbox{if } \frac{4p-1}{12}<r\leq \frac{10p-1}{12}, \\
      6r^2+r-12pr+(6p^2-p) & \mbox{if } r>\frac{10p-1}{12}.
  \end{cases}
    $$
    Define 
    $$
    N=N(p,i)=\max\left(\bigcup_{s=0}^{p-1}\min(i\mathcal{L}(r):\,i\mathcal{L}(r)\equiv s\pmod{p},\,0\leq r\leq p-1)\right)-p.
    $$
  Then for $n>N$,
   \begin{align*}
       a_n\begin{cases}
           >0, &\mbox{if $n\equiv i(6r^2+r) \pmod p$ with $r\leq \frac{4p-1}{12}$ or $r>\frac{10p-1}{12}$,}\\
           <0, &\mbox{if $n\equiv i(6r^2+r) \pmod p$ with $\frac{4p-1}{12} <r\leq \frac{10p-1}{12}$,}\\
           =0, &\mbox{if  $n\neq i(6r^2+r) \pmod p$.}
       \end{cases}
   \end{align*}
  
\end{theorem}

\begin{proof}

First of all, following the notation in Corollary~\ref{qqdissec} we prove that if $6r^2+r \equiv6r'^2+r' \equiv t \pmod p,\text{ then } (-1)^{s(r)}=(-1)^{s(r')} $.
$\text{ Since }p=3k \pm 1, \text{then if }k=2t+1, p=6t+3\pm1$ is even; however, $p$ is prime. So one can just assume $p=6t\pm1$. When $p=6t+1 \text{ and } r,r' \in \mathbb{Z}/p\mathbb{Z}$,
\begin{align*}
    6r^2+r \equiv6r'^2+r' \pmod p \Leftrightarrow (6r+6r'+1)(r-r')\equiv0 \pmod p.
\end{align*}
When $p|(r-r')$, it is done.
If $p  \mid (6r+6r'+1)$, then
\begin{align*}
r+r'\equiv\frac{p-1}{6}\equiv\frac{14p-2}{12} \pmod p. 
\end{align*}
Thus, one can conclude that if $r<\frac{4p-1}{12}$, then $r'>\frac{10p-1}{12}$, and similarly if $ r>\frac{4p-1}{12}$, then $r'<\frac{10p-1}{12}$. So $(-1)^{s(r)}=(-1)^{s(r')} $.
Similarly, when $p=6t-1 \text{ and } r,r' \in \mathbb{Z}/p\mathbb{Z}$, one can find that
\begin{align*}
r+r'\equiv\frac{3t-1}{3}=t-3^{-1}\equiv 5t-1 \pmod {p},
\end{align*}
as $(-3)^{-1}\equiv4t-1 \pmod{p}$, and thus,
 if $r<\frac{2m-1}{12}=\frac{4t-1}{4}$,then $r'>\frac{8m-1}{12}=\frac{16t-3}{4}$ and if $ r>\frac{2m-1}{12}$,then $r'<\frac{8m-1}{12}$. So $(-1)^{s(r)}=(-1)^{s(r')} $.

Recall by Corollary~\ref{qqdissec}, one has that
   $$(q;q)_{\infty}=\sum_{r=0}^{p-1}(-1)^{s(r)} q^{\mathcal{L}(r)}(q^{t_{1}(r)},q^{4p^{2}-t_{1}(r)},q^{4p^{2}};q^{4p^{2}})_{\infty}(q^{t_{2}(r)},q^{8p^{2}-t_{2}(r)};q^{8p^{2}})_{\infty},$$
 where when $p\equiv 1 \pmod 3$, 
\begin{align*}
  t_{1}(r)&=
  \begin{cases}
      \frac{2p^2+p}{3}+4pr, & \mbox{if } r<\frac{10p-1}{12}, \\
      \frac{-10p^2+p}{3}+4pr, & \mbox{if } r>\frac{10p-1}{12}.
  \end{cases} \\
  t_{2}(r)&=
  \begin{cases}
      \frac{16p^2+2p}{3}+8pr, & \mbox{if } r<\frac{4p-1}{12}, \\
      \frac{-8p^2+2p}{3}+8pr, & \mbox{if } r>\frac{4p-1}{12}.
  \end{cases} \\
  \mathcal{L}(r)&=
  \begin{cases}
      6r^2+r    & \mbox{if } r\leq \frac{4p-1}{12}, \\
      6r^2+r-8pr+\frac{8p^2-2p}{3} & \mbox{if } \frac{4p-1}{12}<r\leq \frac{10p-1}{12}, \\
      6r^2+r-12pr+(6p^2-p) & \mbox{if } r>\frac{10p-1}{12}.
  \end{cases}\\
 s(r)&=
 \begin{cases}
   0, & \mbox{if } r\leq \frac{4p-1}{12}, \\
  1, & \mbox{if } \frac{4p-1}{12}<r\leq \frac{10p-1}{12}, \\
  2, & \mbox{if } r>\frac{10p-1}{12}.
 \end{cases}
\end{align*}
Note that if $t_1(r)=t_2(r)$, then one must have 
$$
  \begin{cases}
      \frac{2p^2+p}{3}+4pr=\frac{16p^2+2p}{3}+8pr, & \mbox{if } r<\frac{4p-1}{12}, \\
      \frac{2p^2+p}{3}+4pr=\frac{-8p^2+2p}{3}+8pr, & \mbox{if } \frac{4p-1}{12}<r<\frac{10p-1}{12}, \\
      \frac{-10p^2+p}{3}+4pr=\frac{-8p^2+2p}{3}+8pr, & \mbox{if } r>\frac{10p-1}{12}.
  \end{cases}
$$
It is straightforward to check that, subject to the assumptions on~$r$, all three of the equations have no solutions.
So, $ t_1(r) \neq t_2(r)$. Also, it is routine to check that the quantities $t_1(r),\, 4p^2-t_1(r),\, 4p^2,\,4p^2+t_1(r),\, 8p^2-t_1(r),\, 8p^2,\,t_2(r),\, 8p^2-t_2(r) $ are pairwise distinct and are all multiples of~$p$. 
Therefore, one can find that
   \begin{align*}
       \frac{(q^i;q^i)_{\infty}}{(q^p;q^p)_{\infty}}&=\sum_{r=0}^{p-1}(-1)^{s(r)} q^{i\mathcal{L}(r)}\frac{(q^{it_{1}(r)},q^{i(4p^{2}-t_{1}(r))},q^{i4p^{2}};q^{i4p^{2}})_{\infty}(q^{it_{2}(r)},q^{i(8p^{2}-t_{2}(r))};q^{i8p^{2}})_{\infty}}{(q^p;q^p)_{\infty}} \\
       &=\sum_{r\leq\frac{4p-1}{12}} q^{i\mathcal{L}(r)}\frac{(q^{it_{1}(r)},q^{i(4p^{2}-t_{1}(r))},q^{i4p^{2}};q^{i4p^{2}})_{\infty}(q^{it_{2}(r)},q^{i(8p^{2}-t_{2}(r))};q^{i8p^{2}})_{\infty}}{(q^p;q^p)_{\infty}}\\
       &\quad-\sum_{\frac{4p-1}{12}<r\leq \frac{10p-1}{12}} q^{i\mathcal{L}(r)}\frac{(q^{it_{1}(r)},q^{i(4p^{2}-t_{1}(r))},q^{i4p^{2}};q^{i4p^{2}})_{\infty}(q^{it_{2}(r)},q^{i(8p^{2}-t_{2}(r))};q^{i8p^{2}})_{\infty}}{(q^p;q^p)_{\infty}}\\
       &\quad+\sum_{r>\frac{10p-1}{12}} q^{i\mathcal{L}(r)}\frac{(q^{it_{1}(r)},q^{i(4p^{2}-t_{1}(r))},q^{i4p^{2}};q^{i4p^{2}})_{\infty}(q^{it_{2}(r)},q^{i(8p^{2}-t_{2}(r))};q^{i8p^{2}})_{\infty}}{(q^p;q^p)_{\infty}}.
   \end{align*} 
   Since the infinite products in the numerator have no factor of~$(1-q^{p})$, then by the definition of~$\mathcal{L}(r)$ and Lemma~\ref{technical} one can deduce that $a_{n}=0$ if $n\not\equiv i(6r^{2}+r)\pmod{p}$, and when $n\equiv i(6r^{2}+r)\pmod{p}$ and $n\geq i\mathcal{L}(r)$, 
   $$
    a_n\begin{cases}
           >0, &\mbox{if $n\equiv i(6r^2+r) \pmod p$ with $r\leq \frac{4p-1}{12}$ or $r>\frac{10p-1}{12}$,}\\
           <0, &\mbox{if $n\equiv i(6r^2+r) \pmod p$ with $\frac{4p-1}{12} <r\leq \frac{10p-1}{12}$.}
       \end{cases}
       $$
       It is clear that for $0\leq s\leq p-1$ and $m_{s}:=\min(i\mathcal{L}(r):\,i\mathcal{L}(r)\equiv s\pmod{p},\,0\leq r\leq p-1)$, if the sign-pattern holds for $a_{m_{s}}$, then it holds for any $n$ in the arithmetic progression $m_{s}+kp$ for $k\geq0$. Therefore, for $N=\max(m_{s}:\, 0\leq s\leq p-1)-p$ and $n>N$, all the minimal indices such that the sign pattern holds have been ranged over, and hence, the sign pattern holds for any $n>N$. 

\end{proof}

\section{Proof of Theorem~\ref{main2}}

In this section, we give the proof of Theorem~\ref{main2}. This will be delivered case by case.

\begin{proof}[Proof of Theorem~\ref{main2} (1)]
Start with the identity \eqref{eta2}
   $$
    \frac{(q^2; q^2)_{\infty}^2}{(q; q)_{\infty}}= \sum_{n=0}^{\infty} q^{n(n+1)/2}:=\sum_{n=0}^{\infty}b_{n}q^{n}.
 $$
Notice that for $n\not\equiv \frac{r(r+1)}{2}$ for any $r\in\mathbb{Z}/p\mathbb{Z}$, the coefficient $b_{n}=0$, and $b_{r(r+1)/2}=1$, and thus, by Lemma~\ref{technical} and the same reasoning used at the end of the proof of Theorem~\ref{main1} for
$$
 \frac{(q^2;q^2)^2_{\infty}}
    {(q;q)_{\infty}(q^{p};q^{p})_{\infty}}=\sum_{n=0}^{\infty}a_{n}q^{n},
$$
one can conclude that $a_{n}\geq1$ when $n>N$ and $n\equiv \frac{r(r+1)}{2}\pmod{p}$ for some $r\in\mathbb{Z}/p\mathbb{Z}$.

\end{proof}



      
 
\begin{proof}[Proof of Theorem~\ref{main2} (2)]
Recall by \eqref{eta1} that
$$
\frac{(q;q)_{\infty}^{2}}{(q^{2};q^{2})_{\infty}} = \sum_{n=-\infty}^{\infty} (-1)^n q^{n^2}:=\sum_{n=0}^{\infty}b_{n}q^{n}.
$$
Then it clear that $b_{n}=0$ when $n\not\equiv t^{2}\pmod{4p}$, equivalently, $n\equiv2,3\pmod{4}$ or $\left(\frac{n}{p}\right)=-1$, and 
  \begin{align*}
        b_n=\begin{cases}
                \leq 0&\mbox{if $n\equiv 4t^2+4t+1\pmod{4p}$,}\\
                \geq 0&\mbox{if $n\equiv 4t^2\pmod{4p}$},
            \end{cases}
        \end{align*}
        for some $t\in\mathbb{Z}/p\mathbb{Z}$.
In particular, for $0\leq t\leq p-1$, $b_{(2t+1)^{2}}=-1<0$ and $b_{(2t)^{2}}=1>0$. So, dividing both sides by $(q^{4p};q^{4p})_{\infty}$ and applying Lemma~\ref{technical} and the same reasoning used in the proof of Theorem~\ref{main1}, one can deduce that for $n>N$, 
\begin{align*}
        a_n=\begin{cases}
                < 0&\mbox{if $n\equiv 4t^2+4t+1\pmod{4p}$,}\\
                > 0&\mbox{if $n\equiv 4t^2\pmod{4p}$},
            \end{cases}
        \end{align*}
        for some $t\in\mathbb{Z}/p\mathbb{Z}$.
\end{proof}



 

\begin{proof}[Proof of Theorem~\ref{main2} (3)]
By Theorem~\ref{3dissec} one has
\begin{align*}
    \frac{(q; q)_\infty^3}{(q^{3};q^{3})_{\infty}^{2}} &= \frac{1}{(q^{3};q^{3})_{\infty}}\left(1+6\sum_{n=1}^{\infty}q^{3n}\frac{1-q^{3n}}{1-q^{9n}}\right)-3q \frac{(q^9; q^9)_\infty^3}{(q^{3};q^{3})_{\infty}}\times \frac{1}{(q^{3};q^{3})_{\infty}}\\
    &=\left(\frac{1}{(q^{3};q^{3})_{\infty}}+6\sum_{n=1}^{\infty}q^{3n}\frac{1-q^{3n}}{(1-q^{9n})(q^{3};q^{3})_{\infty}}\right)-3q \frac{(q^9; q^9)_\infty^3}{(q^{3};q^{3})_{\infty}}\times \frac{1}{(q^{3};q^{3})_{\infty}}.
\end{align*}
Clearly, the coefficient $a_{n}=0$ for $n\equiv2\pmod{3}$, and by the first component on the right hand side together with Lemma~\ref{technical}, it can be found that the coefficient $a_{n}\geq1$ for $n\equiv0\pmod{3}$. 
Finally, recall by the beginning of the proof of Lemma~\ref{lem9} that 
$$
\frac{(q^3;q^3)^3_\infty}{(q;q)_\infty}=1+O(q)
$$
has non-negative coefficients, and thus, so does
$$
\frac{(q^9;q^9)^3_\infty}{(q^{3};q^{3})_\infty}=1+O(q^{3}).
$$
This together with Lemma~\ref{technical} implies that $a_{n}<0$ for $n\equiv1\pmod{3}$ and completes the proof.


\end{proof}


    
 



    \begin{proof}[Proof of Theorem~\ref{main2} (4)]
By~\eqref{eta6}, one has
\begin{align*}
    \frac{(q;q)^2_{\infty}(q^6;q^6)_\infty}
{(q^2;q^2)_\infty(q^3;q^3)_{\infty}} &= \sum_{k =0}^{\infty} q^{((6k+1)^2-1)/8} - 2 \sum_{k=0}^{\infty} q^{((6k+3)^2-1)/8} + \sum_{k=0}^{\infty} q^{((6k+5)^2-1)/8} \\
&=\sum_{k=0}^{\infty}q^{\frac{3k(3k+1)}{2}}-2\sum_{k=0}^{\infty}q^{\frac{(3k+1)(3k+2)}{2}}+\sum_{k=0}^{\infty}q^{\frac{(3k+2)(3k+3)}{2}}\\
&=\sum_{n=0}^{\infty}c_{n}q^{3n}-2q\sum_{n=0}^{\infty}b_{n}q^{3n},
\end{align*}
where clearly $c_{n},b_{n}\geq0$ for any~$n$, and $c_{0}=b_{0}=1$. So, divide both sides by $(q^{6};q^{6})_{\infty}(q^{3};q^{3})_{\infty}$ and notice by Lemma~\ref{technical} that if
$$
\frac{1}{(q^{6};q^{6})_{\infty}(q^{3};q^{3})_{\infty}}=\sum_{n=0}^{\infty}d_{n}q^{3n},
$$
then $d_{n}\geq1$. This together with the signs of $c_{n},b_{n}$ and the fact $c_{0}=b_{0}=1$ implies the desired conclusion.
\end{proof}


\begin{proof}[Proof of Theorem~\ref{main2} (5)]

Start with squaring~\eqref{eta1} to get
\begin{align*}
         \frac{(q,q)_\infty^4}{(q^2;q^2)_\infty^2}=\sum_{n=-\infty}^{\infty}\sum_{m=-\infty}^{\infty}(-1)^{n+m}q^{n^{2}+m^2}=\sum_{n=0}^{\infty}b_{n}q^{n}.
    \end{align*}
    Clearly, $b_{n}=0$ when $n\equiv3\pmod{4}$, which cannot be written as a sum of two squares.
   Also notice that if $n_1^2+m_1^2 \equiv n_2^2+m_2^2 \pmod 4, \text{ then }(n_1+n_2)(n_1-n_2)\equiv(m_1+m_2)(m_1-m_2) \pmod 2 $, and thus, $n_1+m_1\equiv n_2+m_2 \pmod 2$. So one can deduce that
   $$
   b_{n}\begin{cases}
       \geq0 &\mbox{if $n\equiv0,2\pmod{4}$,}\\
       \leq 0&\mbox{if $n\equiv1\pmod{4}$.}
   \end{cases}
   $$
   Finally, since $b_{0}=b_{1}=b_{2}=1$, then dividing both sides by $(q^{4};q^{4})_{\infty}$ and invoking Lemma~\ref{technical} yield the desired conclusion.
\end{proof}


\begin{proof}[Proof of Theorem~\ref{main2} (6)]
Squaring the identity~\eqref{eta5}
$$
 \frac{(q^2;q^2)^{5}_{\infty}}{(q;q)_{\infty}^2(q^4;q^4)_{\infty}^2}
    =\sum_{n=-\infty}^{\infty}q^{n^2},
$$
one has that
    $$
    \frac{(q^2;q^2)^{10}_{\infty}}{(q;q)_{\infty}^4(q^4;q^4)_{\infty}^4}
    =\sum_{n=-\infty}^{\infty}\sum_{m=-\infty}^{\infty}q^{m^2+n^2}=\sum_{n=0}^{\infty}b_{n}q^{n}.
    $$
Clearly, since $m^2+n^2\equiv0,1,2\pmod{4}$, then $b_{n}=0$ if $n\equiv3\pmod{4}$, and $b_{n}\geq0$ for $n\equiv0,1,2\pmod{4}$. In particular, $b_{0}=1$, $b_{1}=2$, $b_{2}=1$. So, after dividing both sides by $(q^{4};q^{4})_{\infty}$, the conclusion follows from these together with Lemma~\ref{technical}.
\end{proof}


\begin{proof}[Proof of Theorem~\ref{main2} (7)]
      Recall the following classical identity of Ramanujan \cite{H}:
    \begin{align*}
        (q; q)_{\infty} = (q^{25}; q^{25})_{\infty} \left[ \frac{(q^{10}, q^{15}; q^{25})_{\infty}}{(q^5, q^{20}; q^{25})_{\infty}} - q - q^2 \frac{(q^5, q^{20}; q^{25})_{\infty}}{(q^{10}, q^{15}; q^{25})_{\infty}} \right].
    \end{align*}
    Squaring both sides gives that
    \begin{align*}
        (q; q)^2_{\infty} 
        &= (q^{25}; q^{25})^2_{\infty} \left[ \left( \frac{(q^{10}, q^{15}; q^{25})_{\infty}}{(q^5, q^{20}; q^{25})_{\infty}} \right)^2 - 2q \frac{(q^{10}, q^{15}; q^{25})_{\infty}}{(q^5, q^{20}; q^{25})_{\infty}}-q^{2} \right. \\
        &\qquad \qquad\qquad\qquad+ 2q^3  \frac{(q^5, q^{20}; q^{25})_{\infty}}{(q^{10}, q^{15}; q^{25})_{\infty}} \left. + q^4  \left( \frac{(q^5, q^{20}; q^{25})_{\infty}}{(q^{10}, q^{15}; q^{25})_{\infty}} \right)^2 \right].
    \end{align*}    
    Note by Lemma~\ref{technical} that after dividing both sides by~$(q^{5};q^{5})_{\infty}^{3}$, each of the resulting infinite-product components is of the form
    $$
    \prod_{n\in S}\frac{1}{(1-q^{5n})^{2}}\times\frac{1}{(q^{5};q^{5})_{\infty}}=\left(1+\sum_{n=1}^{\infty}c_{n}q^{5n}\right)\frac{1}{(q^{5};q^{5})_{\infty}}=\sum_{n=0}^{\infty}A_{n}q^{5n},
    $$
    with $c_{n}\geq0$ for any~$n$ in some subset $S$ of $\mathbb{Z}_{>0}$, and therefore, $A_{n}\geq1$ for any~$n$. From this the desired conclusion follows.
\end{proof}


\begin{proof}[Proof of Theorem~\ref{main2} (8)]
Upon Lemma~\ref{lem9}, it suffices to prove that
$$
a_{3n}\begin{cases}
    >0&\mbox{if $n\equiv0\pmod{3}$,}\\
    <0&\mbox{if $n\equiv1\pmod{3}$,}\\
    =0&\mbox{if $n\equiv 2\pmod{3}$.}
\end{cases}
$$
Notice by the proof of Lemma~\ref{lem9} that
\begin{multline*}
    \frac{(q;q)^9}{(q^3;q^3)^9}=\\
    \frac{1}{(q^3;q^3)^{6}}\left(a^3(q^3)-9qa^2(q^3)\frac{(q^{9};q^{9})_{\infty}^{3}}{(q^{3};q^{3})_{\infty}}+27q^{2}a(q^3)\left(\frac{(q^{9};q^{9})_{\infty}^{3}}{(q^{3};q^{3})_{\infty}}\right)^2-27q^3\left(\frac{(q^{9};q^{9})_{\infty}^{3}}{(q^{3};q^{3})_{\infty}}\right)^3\right),
\end{multline*}
and thus,
$$
\sum_{n=0}^{\infty}a_{3n}q^{n}=\frac{a^3(q)}{(q;q)^{6}}-27q\frac{(q^{3};q^{3})_{\infty}^{9}}{(q;q)_{\infty}^{9}}.
$$
Making use of the identities given in~\eqref{abc1} and Theorem~\ref{3dissec}, it is straightforward to deduce that
\begin{align*}
\sum_{n=0}^{\infty}a_{3n}q^{n}&=\frac{(q;q)_{\infty}^{3}}{(q^{3};q^{3})_{\infty}^{3}}\\
&=\frac{1}{(q^3;q^3)_\infty^{2}}\left(1+6\sum_{n=1}^{\infty}q^{3n}\frac{1-q^{3n}}{1-q^{9n}}\right)-3q \frac{(q^9; q^9)_\infty^3}{(q^{3};q^{3})_{\infty}^{3}}.
\end{align*}
Finally, from this together with Lemma~\ref{technical} the desired sign-pattern of period~$3$ for $a_{3n}$
follows.    
\end{proof}


\begin{proof}[Proof of Theorem~\ref{main2} (9)]
By Lemma~\ref{lem9} one already has 
$$
a_{n}\begin{cases}
    >0&\mbox{if $n\equiv2\pmod{3}$, i.e., $n\equiv 2,5,8\pmod{9}$,}\\
    <0&\mbox{if $n\equiv1\pmod{3}$, i.e., $n\equiv 1,4,7\pmod{9}$.}
\end{cases}
$$
It remains to treat the case $n\equiv0\pmod{3}$. To this end, by the proof of Lemma~\ref{lem9} and making use of~\eqref{abc1}, one can first deduce
\begin{align*}
     \frac{(q;q)_{\infty}^9}{(q^3;q^3)_{\infty}^i}
     &=\frac{(q^{3};q^{3})_{\infty}^{12-i}}{(q^{9};q^{9})_{\infty}^{3}}-9q\frac{a(q^{3})^{2}}{(q^3;q^3)_{\infty}^{i-3}}\left(\frac{(q^{9};q^{9})_{\infty}^{3}}{(q^{3};q^{3})_{\infty}}\right)+27q^{2}\frac{a(q^{3})}{(q^3;q^3)_{\infty}^{i-3}}\left(\frac{(q^{9};q^{9})_{\infty}^{3}}{(q^{3};q^{3})_{\infty}}\right)^{2},
\end{align*}
and thus,
$$
\sum_{n=0}^{\infty}a_{3n}q^{3n}=\frac{(q^{3};q^{3})_{\infty}^{12-i}}{(q^{9};q^{9})_{\infty}^{3}}.
$$

When $i>12$, it is clear by Lemma~\ref{technical} that $a_{3n}\geq1$ for any~$n$, and this justifies part (c).

When $i=12$, by Lemma~\ref{technical} one can find that $a_{3n}=0$ for $n\equiv1,2\pmod{3}$ and $a_{3n}>0$ for $n\equiv0\pmod{3}$, i.e., $a_{n}=0$ for $n\equiv3,6\pmod{9}$ and $a_{n}>0$ for $n\equiv0\pmod{9}$. This justifies part (b).

Finally, when $i=11$, first note by Theorem~\ref{3dissec} that
\begin{align*}
     \frac{(q^{3};q^{3})_{\infty}}{(q^{9};q^{9})_{\infty}^{3}}&= \frac{1}{(q^3, q^6, q^9, q^{18}, q^{21}, q^{24}, q^{27}; q^{27})_\infty|_{q\to q^{3}}}\times\frac{1}{(q^{9};q^{9})_{\infty}^{2}}\\
     &\quad- q^{3}\frac{1}{(q^3, q^9, q^{12}, q^{15}, q^{18}, q^{24}, q^{27}; q^{27})_\infty|_{q\to q^{3}}}\times\frac{1}{(q^{9};q^{9})_{\infty}^{2}}\\
    &\quad- q^{6}\frac{1}{(q^6, q^9, q^{12}, q^{15}, q^{18}, q^{21}, q^{27}; q^{27})_\infty|_{q\to q^{3}}}\times\frac{1}{(q^{9};q^{9})_{\infty}^{2}}.
\end{align*}
 So once again, by Lemma~\ref{technical}, one can find that
 $$
 a_{n}\begin{cases}
     >0&\mbox{if $n\equiv0\pmod{9}$,}\\
     <0&\mbox{if $n\equiv3,6\pmod{9}$.}
 \end{cases}
 $$
This completes the proof.

\end{proof}

\section{Concluding Remarks}

In \cite{W22} Wang considered powers of the infinite Borwein product, i.e. eta quotients of the form
\[
G_t^m(q):=
\left(
\frac{(q;q)_{\infty}}{(q^t;q^t)_{\infty}}
\right)^m
=:
\sum_{n=0}^\infty c_t^{(m)}(n)q^n,
\]
for positive integers $m$ and $t$.
Note that $G_t(q)$ was previously defined for $t=p$, $p$ a prime, at \eqref{GPeq1}. 

For  $t$ and $m$  positive integers satisfying $m(t-1)\leq 24$, Wang gave asymptotic formula for $c_t^{(m)}(n)$, and give characterizations of the $n$ for which $c_t^{(m)}(n)$ is positive, negative or zero. It was also shown that $c_t^{(m)}(n)$ is ultimately periodic in sign. Wang also conjectured that for any positive integers $t$ and $m$, that the sequence of coefficients $\{c_t^{(m)}(n)\}_{n\geq 0}$ is ultimately periodic in sign, with least period of sign divisible by $t$. Several cases of this conjecture were proven.

In a subsequent paper we intend to investigate what the methods used in the present paper can contribute to this conjecture and similar conjectures for the related infinite products 
\[
\left(
\frac{(q^j;q^j)_{\infty}}{(q^t;q^t)_{\infty}}
\right)^m, \hspace{25pt}
\frac{(q^j;q^j)_{\infty}^m}{(q^t;q^t)_{\infty}}.
\]
If $\gcd(t,6)=1$, then Corollary \ref{qqdissec} can be used to get the $t$-dissection of $(q^j;q^j)_{\infty}$
and hence of  $(q^j;q^j)_{\infty}^m/(q^t;q^t)_{\infty}$ or $(q^j;q^j)_{\infty}^m/(q^t;q^t)_{\infty}^m$ . Of course any particular term in this $t$-dissection 
will in general not be a single infinite product, and instead will be a sum of infinite products with mixed signs, and thus, in contrast to the situation in the present paper,  not so obvious that coefficients will ultimately be  all of the same sign or ultimately periodic in sign.

As part of a preliminary investigation into this new situation, we looked at the eta quotients $A(q):=(q^2;q^2)_{\infty}^5/(q^7;q^7)_{\infty}$  and $B(q):=(q^3;q^3)_{\infty}^5/(q^7;q^7)_{\infty}$, and counted the number of coefficients (in the first 50,000) that were negative, zero and positive in each of the arithmetic progressions $r \pmod 7$, $r=0,1,\dots 6$. These counts are shown in Table \ref{ta1}.

{\allowdisplaybreaks
\begin{center}
\begin{longtable}{|C|CCC||C|CCC|}
\hline
&&&&&&&\\
&&A(q)&& &&B(q)&\\
&&&&&&&\\
\hline
r& \# - & \# 0 & \# + &r& \# - & \# 0 & \# + \\
\hline
 0 & 0 & 0 & 7142 & 0 & 0 & 0 & 7142 \\
 1 & 7141 & 1 & 0 & 1 & 7140 & 2 & 0 \\
 2 & 3319 & 504 & 3319 & 2 & 0 & 1 & 7141 \\
 3 & 7141 & 1 & 0 & 3 & 3300 & 518 & 3324 \\
 4 & 3285 & 507 & 3350 & 4 & 3294 & 525 & 3323 \\
 5 & 3279 & 509 & 3354 & 5 & 7141 & 1 & 0 \\
 6 & 0 & 0 & 7142 & 6 & 3292 & 524 & 3326 \\
\hline
\caption{Counts of negative, zero and positive coefficients in the arithmetic progressions $7 n+r$, $0\leq r \leq 6$ in the series expansion of $A(q)$ and $B(q)$ }\label{ta1}
\end{longtable}
\end{center}
}

It can be seen that while the evidence suggests that coefficients in some arithmetic progressions are ultimately all of the same sign ($r \equiv 0,1,3,6 \pmod 7$ for $A(q)$ and $r \equiv 0,1,2,5 \pmod 7$ for $B(q)$), the situation is more complex for the other arithmetic progression. 

In these other arithmetic progressions evidence suggests that positive and negative coefficients are roughly equinumerous, but that in addition there are large numbers of zero coefficients, roughly one fourteenth of the total, so that the distribution of the signs of the coefficients is clearly more complicated.

{\bf Declaration of competing interest.}
On behalf of all authors, the corresponding author states that there is no conflict of interest.

{\bf Data availability.}
No data was used for the research described in the article

\end{document}